\documentclass{amsart}
\usepackage{amsmath,amssymb}

\newtheorem{theorem}{Theorem}[section]
\newtheorem{proposition}[theorem]{Proposition}

\newtheorem{corollary}[theorem]{Corollary}

\numberwithin{equation}{section}
\providecommand{\bysame}{\leavevmode\hbox to3em{\hrulefill}\thinspace}

\def\DJ{{\hbox{D\kern-.8em\raise.15ex\hbox{--}\kern.35em}}}
\def\DJo{$\;$\kern-.4em
    \hbox{D\kern-.8em\raise.15ex\hbox{--}\kern.35em okovi\'c}}

\renewcommand{\subjclassname}{\textup{2000} Mathematics Subject
Classification}

\def\al{{\alpha}}

\def\bZ{{\mbox{\bf Z}}}

\begin{document}

\title {On the base sequence conjecture}

\author {Dragomir \v{Z}. \DJo}
\address{Department of Pure Mathematics, University of Waterloo,
Waterloo, Ontario, N2L 3G1, Canada}
\email{djokovic@uwaterloo.ca}

\keywords{Base sequences, normal sequences, near-normal sequences,
T-sequences, Yang numbers}
\date{}


\begin{abstract}
Let $BS(m,n)$ denote the set of base sequences $(A;B;C;D)$,
with $A$ and $B$ of length $m$ and $C$ and $D$ of length $n$.
The base sequence conjecture (BSC) asserts that $BS(n+1,n)$ 
exist (i.e., are non-empty) for all $n$. 
This is known to be true for $n\le36$ and when $n$ is a Golay number. 
We show that it is also true for $n=37$ and $n=38$. 
It is worth pointing out that BSC is stronger than the famous 
Hadamard matrix conjecture.

In order to demonstrate the abundance of base sequences, we have 
previously attached to $BS(n+1,n)$ a graph $\Gamma_n$
and computed the $\Gamma_n$ for $n\le27$. We now extend these 
computations and determine the $\Gamma_n$ for $28\le n\le 35$. 
We also propose a conjecture describing these graphs in general.
\end{abstract}

\maketitle \subjclassname{ 05B20, 05B30 }

\section{Introduction}

By {\em binary} respectively {\em ternary sequence} we mean a 
sequence $A=a_1,a_2,\ldots,a_m$ whose terms belong to $\{\pm1\}$ 
respectively $\{0,\pm1\}$. 
To such a sequence we associate the polynomial 
$A(z)=a_1+a_2z+\cdots+a_mz^{m-1}$.
We refer to the Laurent polynomial $N(A)=A(z)A(z^{-1})$ as 
the {\em norm} of $A$.
{\em Base sequences} $(A;B;C;D)$ are quadruples of binary sequences,
with $A$ and $B$ of length $m$ and $C$ and $D$ of length $n$, 
and such that 
\begin{equation} \label{norm}
N(A)+N(B)+N(C)+N(D)=2(m+n).
\end{equation}
(The last condition is equivalent to the vanishing of the sum of
the aperiodic auto-correlation functions of $A,B,C$ and $D$.)
We denote the set of such base sequences by $BS(m,n)$.
Base sequences, and their special cases such as normal and near-normal
sequences, play an important role in the construction of Hadamard
matrices \cite{HCD,SY}. For instance, the recent discovery of a 
Hadamard matrix of order 428 \cite{KT} used a $BS(71,36)$, 
constructed specially for that purpose.

As explained in \cite{DZ1}, we can view the normal sequences $NS(n)$ 
and near-normal sequences $NN(n)$ as subsets of $BS(n+1,n)$. 
For normal sequences $2n$ must be a sum of three squares, and 
for near-normal sequences $n$ must be even or 1.
The base sequences $(A;B;C;D)\in BS(n+1,n)$ are {\em normal}
respectively {\em near-normal} if $b_i=a_i$
respectively $b_i=(-1)^{i-1}a_i$ for all $i\le n$.

The {\em base sequence conjecture} (BSC), first proposed explicitly
in \cite{DZ1} (see also \cite{HCD}), asserts that the 
$BS(n+1,n)$ exist for all integers $n\ge0$. Implicitly, it appears
in earlier papers of J. Seberry and C.H. Yang. So far, BSC has
been verified for all $n\le36$ and it is also well known that it
holds when $n$ is a {\em Golay number}, i.e., when 
$n=2^a 10^b 26^c$ where $a,b,c$ are nonnegative integers.
For the cases $n\le32$ and references to previous work by other
authors see \cite{DZ1,KKJ,KKS,KJ,SY}. For the cases $n=33,34,35$
see \cite{KS} or the Tables 7-9 below, and for $n=36$ 
see \cite{DZ3} or the next section and the Table 10 below. 

{\em $T$-sequences} are quadruples of ternary sequences, $(X;Y;Z;W)$,
all of the same length $n$ such that for each index $i$ exactly one
of the terms $x_i$, $y_i$, $z_i$, $w_i$ is nonzero, and
$$ N(X)+N(Y)+N(Z)+N(W)=n. $$
We denote by $TS(n)$ the set of all $T$-sequences of length $n$.
The {\em $T$-sequence conjecture} (TSC) asserts that $TS(n)$ exist 
for all integers $n\ge1$.

In section \ref{BSC} we show that BSC is also valid for 
$n=37$ and $n=38$. Our example for $n=38$ consists of near-normal
sequences. Consequently, the number 77 is a Yang number. We recall 
that {\em Yang numbers} are odd integers $2s+1$ for which $NS(s)$ or
$NN(s)$ is not empty. We also update the status of the TSC.

Let $\al=(A;B;C;D)\in BS(m,n)$ and let $a=A(1),b=B(1),c=C(1),d=D(1)$
and $a^*=A(-1),b^*=B(-1),c^*=C(-1),d^*=D(-1)$. By setting $z=1$
in the norm identity (\ref{norm}), we see that the squares 
$a^2,b^2,c^2,d^2$, 
arranged in decreasing order, form a partition of $2(m+n)$. 
The same is true for the squares of $a^*,b^*,c^*,d^*$.
We denote the former partition by $p_\al$ and the latter by $p^*_\al$.

In the early searches for base sequences $BS(n+1,n)$ 
the objective apparently was to construct, 
for each partition $p$ of $2(2n+1)$ into four squares,
base sequences $\al=(A;B;C;D)\in BS(n+1,n)$ such that $p_\al=p$.
For this we refer the reader to the paper \cite{KKS} and its references. A more ambitious program to construct  
base sequences $\al=(A;B;C;D)\in BS(n+1,n)$ with $p_\al$
and $p^*_\al$ specified was initiated in our paper \cite{DZ1}.
For that purpose, we have defined there the graphs $\Gamma_n$, $n\ge0$. 

In section \ref{Gama} we recall the definition of the $\Gamma_n$. 
These are undirected graphs with loops allowed but no multiple edges. 
They were determined for $n\le27$ by means of extensive computations
of base sequences. We extend these computations to cover the
cases $n=28,29,\ldots,35$. The base sequences that we need are listed
in Tables 2-9 in section \ref{bz}. On the basis of these computations,
we propose a conjecture about the isomorphism types of the graphs
$\Gamma_n$ and show that the conjecture is valid for $n\le35$.

In section \ref{Algoritam} we describe briefly our algorithm for
exhaustive search of the base sequences $BS(n+1,n)$.

In section \ref{NSiNN} we report the results of our recent
searches for $NS(n)$ and $NN(n)$. We also describe what is currently
known about the existence of Yang numbers.

\section{Current status of BSC} \label{BSC}

At the time when BSC was formulated in \cite{DZ1}, it was known that
it holds for $n\le32$. This was extended to $n\le35$ by Kounias
and Sotirakoglou \cite{KS}. The examples of $BSC(n+1,n)$ for $n=36$ 
and $n=38$ were constructed in the course of our exhaustive searches 
for near-normal sequences \cite{DZ3,DZ4}. We have recently constructed
an example for $n=37$. These three examples will be given below.

\begin{proposition} \label{BSC hip}
The base sequences $BS(n+1,n)$ exist for $n\le 38$ (and for
all Golay numbers $n$).
\end{proposition}

As explained above, it suffices to give examples
of $BS(n+1,n)$ for $n=36,37$ and $38$. To avoid possible errors,
we shall give all base sequences in encoded compact form which
is used in our computer program. Although this encoding scheme
has been described in several of our previous papers, we shall give
the details once again for the convenience of the reader.

Let $(A;B;C;D) \in BS(n+1,n)$. We encode the pairs $(A;B)$
and $(C;D)$ separately by using the same scheme.
We decompose the pair $(A;B)$ into quads
$$ \left[ \begin{array}{ll} a_i & a_{n+2-i} \\ 
b_i & b_{n+2-i} \end{array} \right],\quad i=1,2,\ldots,
\left[ \frac{n+1}{2} \right], $$ 
and, if $n=2m$ is even, the central column
$ \left[ \begin{array}{l} a_{m+1} \\ b_{m+1} \end{array} \right]. $
We can assume (and we do) that the first quad of $(A;B)$ is 
$ \left[ \begin{array}{ll} + & + \\ + & - \end{array} \right]. $
We attach to this particular quad the label 0.
The other quads in $(A;B)$ and all the quads of the pair $(C;D)$,
shown with their labels, must be one of the following:
\begin{center}
\begin{eqnarray*}
1=\left[ \begin{array}{ll} + & + \\ + & + \end{array} \right],\quad 
2=\left[ \begin{array}{ll} + & + \\ - & - \end{array} \right],\quad 
3=\left[ \begin{array}{ll} - & + \\ - & + \end{array} \right],\quad 
4=\left[ \begin{array}{ll} + & - \\ - & + \end{array} \right], \\
5=\left[ \begin{array}{ll} - & + \\ + & - \end{array} \right],\quad 
6=\left[ \begin{array}{ll} + & - \\ + & - \end{array} \right],\quad 
7=\left[ \begin{array}{ll} - & - \\ + & + \end{array} \right],\quad 
8=\left[ \begin{array}{ll} - & - \\ - & - \end{array} \right].
\end{eqnarray*}
\end{center}
The central column (if present) is encoded as
$$
0=\left[ \begin{array}{l} + \\ + \end{array} \right],\quad
1=\left[ \begin{array}{l} + \\ - \end{array} \right],\quad
2=\left[ \begin{array}{l} - \\ + \end{array} \right],\quad
3=\left[ \begin{array}{l} - \\ - \end{array} \right].
$$

If $n=2m$ is even, the pair $(A;B)$ is encoded as the
sequence $q_1q_2\ldots q_mq_{m+1}$, where $q_i$,
$1\le i\le m$, is the label of the $i$th quad and $q_{m+1}$
is the label of the central column. If $n=2m-1$ is odd,
then $(A;B)$ is encoded by $q_1q_2\ldots q_m$, where $q_i$
is the label of the  $i$th quad for each $i$.
We use the same recipe to encode the pair $(C;D)$.

As an example, the base sequences
\begin{eqnarray*}
A &=& +,+,+,+,-,-,+,-,+ \\
B &=& +,+,+,-,+,+,+,-,- \\
C &=& +,+,-,-,+,-,-,+ \\
D &=& +,+,+,+,-,+,-,+
\end{eqnarray*}
are encoded as $06142;\, 1675$.

With this notation, the three promised base sequences $BS(n+1,n)$ are
\begin{eqnarray*}
&& n=36: 0764841234846532153; 165154775335162126 \\
&& n=37: 0686287846153524326; 1153175814738523732 \\
&& n=38: 07641237828515856281; 1782612553714317675 
\end{eqnarray*}
Those for $n=36,38$ are in fact near-normal.

It is well known that there exist maps $BS(m,n) \to TS(m+n)$ and 
$TS(n) \to TS(2n)$. By using the Proposition \ref{BSC hip} and taking 
into account the  \cite[Remark V.8.47]{HCD}, we obtain

\begin{corollary}
Apart from the two undecided cases $n=79,97$, the $T$-sequences
$TS(n)$ exist for all $n\le100$.
\end{corollary}

\section{The $\Gamma$-conjecture} \label{Gama}

We begin by recalling the definition of the graph $\Gamma_n$.
Its vertex set is the set of all partitions of $4n+2$ 
into four squares (including 0 and with repetitions allowed).
We postulate that $\Gamma_n$ may have loops but we do not permit
multiple edges. There is a loop at a vertex $p$ if and only if
there exist base sequences $\al\in BS(n+1,n)$ such that
$p_\al=p_\al^*=p$. If $p$ and $q$ are two distinct vertices,
then $\{p,q\}$ is an edge of $\Gamma_n$ if and only if there
exist base sequences $\beta\in BS(n+1,n)$ such that
$\{p_\beta,p_\beta^*\}=\{p,q\}$. This completes the definition of
$\Gamma_n$. We refer to any $\al\in BS(n+1,n)$ as a witness for
the edge $\{p_\al,p_\al^*\}$ of $\Gamma_n$.

While BSC simply asserts that each $BS(n+1,n)$ is non-empty, 
we shall propose a new conjecture which gives the description of
the graphs $\Gamma_n$. 

To state this new conjecture, we need some more notation. Let
$\al$ be as above and assume that $n$ is fixed. 
Note that $a\equiv b\equiv n+1 \pmod{2}$ and
$c\equiv d\equiv n \pmod{2}$. Thus exactly two of the integers
$a,b,c,d$ are even. If $n$ is even, one can show (see \cite{DZ1}) 
that these two even integers are congruent to each other modulo 4. 
In that case we say that the vertex $\al$ is {\em even} respectively 
{\em odd} if they are congruent to 0 respectively 2 modulo 4. 
Thus, for even $n$, the vertex set is partitioned into even and odd vertices.

Let $\nu$ denote the number of vertices of $\Gamma_n$. If $n$ 
is even, let $\nu_0$ respectively $\nu_1$ denote the number of even
respectively odd vertices of $\Gamma_n$. Of course, we have 
$\nu_0+\nu_1=\nu$ when $n$ is even. In Table 1 we give, for $0\le n\le 40$, 
the value of $\nu$ for odd $n$ and the values of $\nu_0$ and $\nu_1$
for even $n$. 

\begin{center}
\begin{tabular}{rrl|rl|rrl|rr}
\multicolumn{10}{c}{Table 1: The numbers $\nu_0$, $\nu_1$ and $\nu$} \\ \hline 
$n$&$\nu_0$&$\nu_1\quad$ & $\quad n$&$\nu\quad$ & 
$\quad n$&$\nu_0$&$\nu_1\quad$ & $\quad n$&$\nu$ \\ \hline
0 & 1 & 0 & 1 & 1 & 20 & 5 & 2 & 21 & 5 \\
2 & 1 & 1 & 3 & 1 & 22 & 5 & 4 & 23 & 4 \\
4 & 2 & 1 & 5 & 2 & 24 & 4 & 3 & 25 & 7 \\
6 & 2 & 1 & 7 & 2 & 26 & 4 & 3 & 27 & 6 \\
8 & 3 & 1 & 9 & 3 & 28 & 5 & 3 & 29 & 7 \\
10 & 2 & 2 & 11 & 2 & 30 & 4 & 4 & 31 & 8 \\
12 & 4 & 1 & 13 & 5 & 32 & 6 & 4 & 33 & 7 \\
14 & 2 & 3 & 15 & 3 & 34 & 5 & 5 & 35 & 5 \\
16 & 4 & 2 & 17 & 5 & 36 & 6 & 3 & 37 & 12 \\
18 & 3 & 2 & 19 & 4 & 38 & 5 & 6 & 39 & 6 \\
   &   &   &    &   & 40 & 9 & 4 &    & \\
\hline
\end{tabular} \\
\end{center}

Let $K_m$ denote the complete graph on $m$ vertices. Any two
distinct vertices are joined by a single edge. However, $K_m$ has
no loops. If we enlarge $K_m$ by attaching a loop at each vertex,
we obtain the graph $K_m^0$. By $K_{m,n}$ we denote the complete
bipartite graph with $m$ respectively $n$ vertices in the first 
respectively second part. 
The disjoint union of two graphs will be written as a sum. \\

{\bf $\Gamma$-conjecture.} {\em $\Gamma_n$ is isomorphic to

(a) $K_\nu^0$ if $n$ is odd;

(b) $K_{\nu_0,\nu_1}$ if $n\equiv2 \pmod{4}$;

(c) $K_{\nu_0}^0+K_{\nu_1}^0$ if $n\equiv0 \pmod{4}$ 
except for $n=4,8,12$. } \\

The graphs $\Gamma_n$ for $n=4,8,12$ are described in \cite{DZ1}.
Since we always have $\nu\ge1$, BSC is a consequence of the
$\Gamma$-conjecture if $n\not\equiv 2 \pmod{4}$. This would also be 
true when $n\equiv 2 \pmod{4}$ provided that one can show that
both $\nu_0$ and $\nu_1$ are nonzero. We can formulate this as 
the following number-theoretical question.

{\bf Question} Let $S=\{k^2:k\in\bZ\}$ respectively $T=\{k(k+1)/2:k\in\bZ\}$
be the set of squares respectively triangular numbers. Let
$S_2=\{x+y:x,y\in S\}$ and $T_2=\{x+y:x,y\in T\}$. Does
the set $\{4x+y:x,y\in T_2\}$ respectively $\{2x+y:x\in S_2,y\in T_2\}$ 
contain all even respectively odd nonnegative integers?

(The BSC implies that the answer is affirmative in both cases.)

We give now the current status of the $\Gamma$-conjecture.

\begin{proposition} \label{Gama hip}
The $\Gamma$-conjecture is valid for $n\le35$.
\end{proposition}

We have to construct witnesses of all hypothetical edges of
$\Gamma_n$. This was accomplished in \cite{DZ1} for $n\le27$, 
while for $n=28$ two witnesses were missing.
Tables 2-9 of the appendix confirm the $\Gamma$-conjecture for 
$n=28,29,\ldots,35$ as they contain witnesses for all 
hypothetical edges of $\Gamma_n$.

We have partial results for $n=36$. Hypothetically, 
$\Gamma_{36}$ has 27 edges. We list the witnesses for 19 of them 
in Table 10.

If $n$ is odd, we use the (decreasing) lexicographic order of 
partitions to enumerate the vertices of $\Gamma_n$. 
If $n$ is even, we enumerate first the even and then the odd 
vertices and arrange them (separately) in the lexicographic 
order. If $n\equiv2 \pmod{4}$ then $\Gamma_n$ is bipartite 
(and there are no loops). The symbol $i$-$j$ in the first column 
of the tables below denotes the edge joining the $i$th and the 
$j$th vertex. If $i=j$, it refers to the loop at the $i$th 
vertex.

For instance, if $n=28$ then there are eight vertices:
\begin{center}
\begin{tabular}{rlcrlcrl}
1) & $(9^2,4^2,4^2,1)$ & \quad & 2) & $(8^2,7^2,1,0)$ & \quad & 3) 
& $(8^2,5^2,5^2,0)$ \\
4) & $(8^2,5^2,4^2,3^2)$ & \quad & 5) & $(7^2,7^2,4^2,0)$ & \quad & & \\
6) & $(10^2,3^2,2^2,1)$ & \quad & 7) & $(9^2,5^2,2^2,2^2)$ & \quad & 8) 
& $(7^2,6^2,5^2,2^2)$ \\
\end{tabular} \\
\end{center}
Since 8,4,0 are all $\equiv0 \pmod{4}$, the first five vertices are even.
Since 10,6,2 are all $\equiv2 \pmod{4}$, the remaining three vertices are odd.
The graph $\Gamma_{28}$ is a disjoint union of $K_5^0$ on even
vertices and $K_3^0$ on odd ones. The first fifteen base sequences
in Table 2 are witnesses for the edges of the ``even'' component 
$K_5^0$, 
and the next six are witnesses for the edges of the ``odd'' component $K_3^0$.

For a witness $\al\in BS(n+1,n)$, the integers $a,b,c,d$ determine 
the vertex $p_\al$ as the partition of $4n+2$ with parts
$a^2,b^2,c^2,d^2$. Similarly, $a^*,b^*,c^*,d^*$ determine the
vertex $p_\al^*$.

\section{Sketch of the algorithm} 
\label{Algoritam}

Our computer program is designed for exhaustive search 
of base sequences $BS(n+1,n)$ for $n\ge7$. The search is divided
into 18 cases by fixing the first three quads of the pair 
$(A;B)$ and the first two quads of $(C;D)$. 
The choice of these cases depends on the parity of $n$.

\begin{center}
\begin{tabular}{llllll}
\multicolumn{6}{c}{Cases for $n$ odd} \\ \hline 
1) 065; 11 & 2) 066; 11 & 3) 068; 11 &
4) 061; 12 & 5) 063; 12 & 6) 064; 12 \\
7) 061; 16 & 8) 063; 16 & 9) 064; 16 &
10) 016; 61 & 11) 017; 61 & 12) 018; 61 \\
13) 016; 64 & 14) 017; 64 & 15) 018; 64 &
16) 011; 66 & 17) 012; 66 & 18) 013; 66 \\
\end{tabular}
\end{center}

\begin{center}
\begin{tabular}{llllll}
\multicolumn{6}{c}{Cases for $n$ even} \\ \hline 
1) 076; 12 & 2) 077; 12 & 3) 078; 12 &
4) 076; 16 & 5) 077; 16 & 6) 078; 16 \\
7) 071; 18 & 8) 072; 18 & 9) 073; 18 &
10) 065; 11 & 11) 066; 11 & 12) 068; 11 \\
13) 061; 12 & 14) 063; 12 & 15) 064; 12 &
16) 061; 16 & 17) 063; 16 & 18) 064; 16 \\
\end{tabular}
\end{center}

Each of the 18 cases is treated separately.
The first quad of the pair $(A;B)$ is always $0$. Thus the 
$n$-th auto-correlation of $(A;B;C;D)$ is 0. The other four 
starting quads are chosen so that the $(n-1)$-st and $(n-2)$-nd 
auto-correlation is $0$. We proceed by selecting the $4$-th quad 
of $(A;B)$ and the $3$-rd quad of $(C;D)$ so that the $(n-3)$-rd 
auto-correlation vanishes. We continue this procedure as far as 
possible. If no selection is possible, we backtrack. If we 
succeed in finding all the quads and the central column, then we 
test whether all the remaining auto-correlations vanish. If not, 
we backtrack. Otherwise we record the base sequences that we 
found. Note that this algorithm does not use any information 
about the possible sums $a,b,c,d$ of the four constituent 
sequences. Thus we do not know in advance what these sums will 
turn out to be.

In order to handle the large values of $n$, say $n>31$, we 
modify the program by breaking it into two phases. The first 
(easy) phase is to collect into a file the initial segments of 
quads, say of length 8 for $(A;B)$ and length 7 for $(C;D)$. 
Such a file has several milions of rows (subcases). It takes 
only several minutes to generate this file. In the second phase 
we use a random number generator to select a row in this file as 
the entry point for our program. The program then completes the 
computation for a fixed number, say $r$, of consecutive rows 
starting from the chosen entry point. We may repeat this 
subroutine, say $s$ times. In our runs, the product $rs$ was 
either 10000 or just 1000. Usually we do not run the program to 
completion as this would require a prohibitively long time. We 
collect all base sequences that the program finds, and stop it 
after 5-6 days. If necessary, we repeat this process several 
times, using different cases, until we find the witnesses for 
all edges of $\Gamma_n$. 

As an example, we mention that the construction of Table 6 took 
in total about 1423 days of CPU time. For this table, we ran the 
parallelized version of our program on two machines at the same time, one used 128 processors at 3.0 GHz and the other 64 
processors at 2.2 GHz. The program constructed in total 2640 
different base sequences $BS(33,32)$.

\section{Recent results on normal and near-normal sequences} 
\label{NSiNN}

We give here a brief summary of our recent results on these two types
of sequences and on Yang numbers.
Let us begin by quoting Theorem V.8.38 from the recent handbook \cite{HCD}.

\begin{theorem} There is no $NS(n)$ for 
$n=6,14,17,21,22,23,24,27,28,30$
(all other orders of $n<31$ exist). $NS(31)$ is the first unknown case.
\end{theorem}

We have carried out exhaustive searches for $NS(n)$ for 
$n=31,33,34,35,36,37,38,39$
and did not find any such sequences. As $32$ and $40$ are Golay 
numbers, we therefore have the following improvement.

\begin{proposition} For $n\le40$, $NS(n)=\emptyset$ if and only if
$$ n \in \{ 6,14,17,21,22,23,24,27,28,30,31,33,34,35,36,37,38,39 \}. $$
The first unknown case is $n=41$.
\end{proposition}

Yang conjecture (see \cite[Conjecture V.8.39]{HCD}) asserts that 
$NN(n)$ exist for all even integers $n$. This has been known to 
be true when $n\le30$ (see \cite{DZ1} and 
\cite[Remark V.8.40] {HCD}). Complete classification of
near-normal sequences has been carried out recently in our notes
\cite{DZ2,DZ3,DZ4} for all even $n\le40$. It turns out that they
exist for all even $n\le40$. Thus Yang conjecture
remains open.

Consequently, we have the following result about Yang numbers
(compare with \cite[Theorem V.8.42.1]{HCD}).

\begin{proposition} For odd integers $n\le81$, $n$ is a Yang number
if and only if
$$ n \not\in \{ 35,43,47,55,63,67,71,75,79 \}. $$
The first unknown case is $n=83$.
\end{proposition}

\section{Acknowledgments}

The author is grateful to NSERC for the continuing support of
his research. Most of this work was made possible by the facilities 
of the Shared Hierarchical Academic Research Computing Network 
(SHARCNET:www.sharcnet.ca). The author also thanks the referees
for their valuable suggestions.

\newpage

\section{Appendix: Lists of base sequences} \label{bz}

\begin{center}
\begin{tabular}{|c|l|l|l|}
\multicolumn{4}{c}{Table 2: $BS(29,28)$} \\ \hline 
Edge & $A$ \& $B$; $C$ \& $D$ & $a,b,c,d$ & $a^*,b^*,c^*,d^*$ \\ \hline
1-1 & $076413275222630$; & $9,-1,4,-4$ & $9,-1,4,4$ \\
& $12875373652226$ &&\\
1-2 & $076514146435673$; & $1,7,8,0$ & $9,-1,4,4$ \\
& $12566715632821$ &&\\
1-3 & $076412161284762$; & $5,5,0,8$ & $1,9,4,4$ \\
& $12876155137475$ &&\\
1-4 & $078482447637733$; & $-9,1,4,-4$ & $-5,-3,4,8$ \\
& $12858753246321$ &&\\
1-5 & $078451311636611$; & $7,7,4,0$ & $-1,-9,4,4$ \\
& $12838752334113$ &&\\
2-2 & $078461443688572$; & $-7,1,0,-8$ & $1,-7,0,8$ \\
& $12848552856354$ &&\\
2-3 & $078457641147620$; & $1,7,0,8$ & $5,-5,8,0$ \\
& $12856747141347$ &&\\
2-4 & $051782353215153$; & $7,1,0,8$ & $3,5,8,4$ \\
& $17678365277211$ &&\\
2-5 & $078485628682111$; & $1,-7,0,-8$ & $-7,-7,0,4$ \\
& $12845558724283$ &&\\
3-3 & $077658617271583$; & $-5,5,8,0$ & $-5,5,0,-8$ \\
& $12852541333416$ &&\\
3-4 & $078466512613430$; & $5,3,4,-8$ & $5,-5,0,8$ \\
& $12862352528373$ &&\\
3-5 & $078517356737323$; & $-5,5,0,8$ & $7,-7,0,4$ \\
& $12747162866717$ &&\\
4-4 & $078458231755712$; & $-3,5,8,-4$ & $5,-3,4,8$ \\
& $12835732236261$ &&\\
4-5 & $078475657853170$; & $-7,7,0,-4$ & $5,3,4,8$ \\
& $12876165548382$ &&\\
5-5 & $078321422423580$; & $7,-7,-4,0$ & $7,-7,4,0$ \\
& $12887533734554$ &&\\
6-6 & $078582621567150$; & $3,1,10,-2$ & $3,1,10,-2$ \\
& $12456332286115$ &&\\
6-7 & $078467557578650$; & $-9,5,2,-2$ & $3,1,10,2$ \\
& $12836515766382$ &&\\
6-8 & $078416634842140$; & $3,1,-2,-10$ & $7,5,2,6$ \\
& $12882758538342$ &&\\
7-7 & $076443181762112$; & $5,9,2,2$ & $9,5,-2,-2$ \\
& $12868357554116$ &&\\
7-8 & $076411216766222$; & $9,5,2,2$ & $5,-7,2,6$ \\
& $12875653427313$ &&\\
8-8 & $078436621518110$; & $7,5,6,-2$ & $7,5,6,2$ \\
& $12886731231325$ &&\\
\hline
\end{tabular} \\
\end{center}

\begin{center}
\begin{tabular}{|c|l|l|l|}
\multicolumn{4}{c}{Table 3: $BS(30,29)$} \\ \hline 
Edge & $A$ \& $B$; $C$ \& $D$ & $a,b,c,d$ & $a^*,b^*,c^*,d^*$ \\ \hline
1-1 & $068362252723438$; & $4,-10,1,-1$ & $10,4,1,-1$ \\
& $118624666538452$ &&\\
1-2 & $068385638777645$; & $-10,0,3,3$ & $4,10,-1,-1$ \\
& $118722343573530$ &&\\
1-3 & $066247531158121$; & $10,4,1,1$ & $0,6,9,1$ \\
& $117543585724280$ &&\\
1-4 & $066217723624145$; & $8,2,7,1$ & $-10,4,-1,1$ \\
& $117432416826461$ &&\\
1-5 & $066417145712627$; & $6,8,3,3$ & $-4,10,-1,-1$ \\
& $117653654785220$ &&\\
1-6 & $068385545252336$; & $2,-8,5,5$ & $-4,10,1,1$ \\
& $118567425535130$ &&\\
1-7 & $066427711368186$; & $2,4,7,7$ & $4,10,-1,-1$ \\
& $117726654641520$ &&\\
2-2 & $066227632114544$; & $10,0,-3,3$ & $0,10,-3,3$ \\
& $117768627585431$ &&\\
2-3 & $066244127461835$; & $6,0,-1,9$ & $0,10,3,-3$ \\
& $117687675413252$ &&\\
2-4 & $066221154863181$; & $10,0,-3,3$ & $8,2,1,7$ \\
& $117586785628152$ &&\\
2-5 & $068246422128374$; & $6,-8,3,3$ & $0,-10,3,3$ \\
& $118657526217580$ &&\\
2-6 & $068248487512863$; & $-2,-8,5,-5$ & $0,-10,-3,3$ \\
& $118768327622521$ &&\\
2-7 & $066325474783574$; & $-4,2,-7,7$ & $10,0,-3,3$ \\
& $117876865367552$ &&\\
3-3 & $066425872412617$; & $6,0,1,9$ & $0,6,1,9$ \\
& $117661785545180$ &&\\
3-4 & $066225637518271$; & $6,0,-1,9$ & $8,2,7,1$ \\
& $117654433817272$ &&\\
3-5 & $066357474847817$; & $-8,6,3,3$ & $6,0,-9,-1$ \\
& $117581625334633$ &&\\
3-6 & $066213624581187$; & $6,0,9,-1$ & $8,2,5,-5$ \\
& $117683252526141$ &&\\
3-7 & $066416178423476$; & $2,4,7,7$ & $0,6,-9,-1$ \\
& $117671223663613$ &&\\
4-4 & $066415721365525$; & $8,2,7,1$ & $2,8,7,1$ \\
& $117726281652351$ &&\\
4-5 & $066417218511536$; & $8,6,3,3$ & $2,8,-1,7$ \\
& $117644754226580$ &&\\
4-6 & $066242378263856$; & $2,-8,7,1$ & $8,2,-5,5$ \\
& $117685486122122$ &&\\
4-7 & $066415277854231$; & $4,2,7,7$ & $-2,8,-1,7$ \\
& $117662161546363$ &&\\
\hline
\end{tabular} \\
\end{center}

\begin{center}
\begin{tabular}{|c|l|l|l|}
\multicolumn{4}{c}{Table 3: (continued)} \\ \hline 
Edge & $A$ \& $B$; $C$ \& $D$ & $a,b,c,d$ & $a^*,b^*,c^*,d^*$ \\ \hline
5-5 & $016186616313366$; & $8,6,3,3$ & $6,8,3,3$ \\
& $641515851514853$ &&\\
5-6 & $066424271211847$; & $8,2,5,5$ & $-6,8,-3,-3$ \\
& $117681267525360$ &&\\
5-7 & $066227415141467$; & $8,6,3,3$ & $2,-4,7,7$ \\
& $117628153854530$ &&\\
6-6 & $068427113134776$; & $2,8,-5,5$ & $-8,2,-5,5$ \\
& $118653736872672$ &&\\
6-7 & $066425635118187$; & $4,2,-7,7$ & $2,8,5,-5$ \\
& $117765384785371$ &&\\
7-7 & $066347444712723$; & $2,4,7,7$ & $4,-2,7,7$ \\
& $117823654415150$ &&\\
\hline
\end{tabular} \\
\end{center}

\begin{center}
\begin{tabular}{|c|l|l|l|}
\multicolumn{4}{c}{Table 4: $BS(31,30)$} \\ \hline 
Edge & $A$ \& $B$; $C$ \& $D$ & $a,b,c,d$ & $a^*,b^*,c^*,d^*$ \\ \hline
1-5 & $0784614381231342$; & $3,3,10,-2$ & $1,-11,0,0$ \\
& $128685615224114$ &&\\
1-6 & $0776853138438782$; & $-11,1,0,0$ & $-1,-9,6,-2$ \\
& $128665371865672$ &&\\
1-7 & $0784216352512611$; & $11,-1,0,0$ & $1,-7,6,6$ \\
& $128863554766615$ &&\\
1-8 & $0784864477847431$; & $-11,1,0,0$ & $-5,-5,6,6$ \\
& $128574476353272$ &&\\
2-5 & $0776162345126151$; & $9,5,0,-4$ & $3,3,10,2$ \\
& $128868657542531$ &&\\
2-6 & $0778853587261780$; & $-9,1,6,2$ & $5,-9,0,4$ \\
& $128558541366151$ &&\\
2-7 & $0778511521651532$; & $5,9,0,-4$ & $7,-1,6,6$ \\
& $128588623471636$ &&\\
2-8 & $0784216213317131$; & $9,5,-4,0$ & $-5,-5,6,6$ \\
& $128863657667445$ &&\\
3-5 & $0764411241717863$; & $3,9,-4,4$ & $-3,3,10,2$ \\
& $128763613567478$ &&\\
3-6 & $0776261117545653$; & $3,9,4,-4$ & $9,-1,6,2$ \\
& $128813253753652$ &&\\
3-7 & $0784162254551610$; & $9,3,4,-4$ & $7,1,6,6$ \\
& $128865236166725$ &&\\
3-8 & $0778565314743723$; & $-5,5,6,6$ & $9,3,4,4$ \\
& $128563665117166$ &&\\
4-5 & $0778586368314251$; & $-3,-3,10,2$ & $3,7,8,0$ \\
& $128566641315214$ &&\\
4-6 & $0512656235371531$; & $9,1,2,6$ & $3,7,0,8$ \\
& $165711846213678$ &&\\
4-7 & $0764323438577832$; & $-7,1,-6,6$ & $7,3,8,0$ \\
& $128767756347465$ &&\\
4-8 & $0564376515151581$; & $5,5,6,6,$ & $7,3,8,0$ \\
& $118772615545132$ &&\\
\hline
\end{tabular} \\
\end{center}

\begin{center}
\begin{tabular}{|c|l|l|l|}
\multicolumn{4}{c}{Table 5: $BS(32,31)$} \\ \hline 
Edge & $A$ \& $B$; $C$ \& $D$ & $a,b,c,d$ & $a^*,b^*,c^*,d^*$ \\ \hline
1-1 & $0653276646881415$; & $2,0,1,11$ & $-2,0,-1,-11$ \\
& $1187615124567762$ &&\\
1-2 & $0653477313582724$; & $0,2,1,11$ & $0,10,-1,5$ \\
& $1186645576741711$ &&\\
1-3 & $0664286361577533$; & $0,2,1,11$ & $4,10,-1,-3$ \\
& $1177645461752362$ &&\\
1-4 & $0664283673814787$; & $-6,0,3,9$ & $-2,0,1,11$ \\
& $1176554618357311$ &&\\
1-5 & $0664483763412781$; & $0,2,1,11$ & $4,2,-5,9$ \\
& $1176755458517611$ &&\\
1-6 & $0663151725347817$; & $2,8,7,-3$ & $2,0,1,11$ \\
& $1178525218466222$ &&\\
1-7 & $0664457618863416$; & $0,2,1,11$ & $-8,-6,-5,1$ \\
& $1176713544647422$ &&\\
1-8 & $0664286477134572$; & $0,2,-1,11$ & $4,6,5,-7$ \\
& $1177653127576363$ &&\\
2-2 & $0663151774538174$; & $0,10,-5,1$ & $0,10,5,-1$ \\
& $1178566324188782$ &&\\
2-3 & $0664134763185177$; & $0,10,5,1$ & $4,10,-1,3$ \\
& $1176835434253213$ &&\\
2-4 & $0663684725887517$; & $-6,0,3,9$ & $10,0,1,-5$ \\
& $1177658146162731$ &&\\
2-5 & $0663554568171527$; & $2,4,-5,9$ & $10,0,5,-1$ \\
& $1177877658254711$ &&\\
2-6 & $0663284811641421$; & $10,0,1,5$ & $2,8,3,7$ \\
& $1178625647454413$ &&\\
2-7 & $0653257763411145$; & $6,8,5,1$ & $10,0,-5,-1$ \\
& $1187716272282540$ &&\\
2-8 & $0664463374577363$; & $-4,6,-5,7$ & $0,10,1,5$ \\
& $1177365857427533$ &&\\
3-3 & $0653271351241777$; & $4,10,3,-1$ & $-4,10,-3,1$ \\
& $1186637254782520$ &&\\
3-4 & $0653485354761371$; & $0,6,9,3$ & $-4,10,3,1$ \\
& $1186373522521312$ &&\\
3-5 & $0663174726216214$; & $10,4,1,3$ & $2,4,-5,9$ \\
& $1178327566525642$ &&\\
3-6 & $0664256357162313$; & $8,2,3,7$ & $4,10,-3,1$ \\
& $1176615635414833$ &&\\
3-7 & $0653461761515422$; & $10,4,1,3$ & $-6,8,-5,1$ \\
& $1187654414627381$ &&\\
3-8 & $0653182153651377$; & $4,6,7,5$ & $4,10,1,3$ \\
& $1186554317646211$ &&\\
\hline
\end{tabular} \\
\end{center}

\begin{center}
\begin{tabular}{|c|l|l|l|}
\multicolumn{4}{c}{Table 5: (continued)} \\ \hline 
Edge & $A$ \& $B$; $C$ \& $D$ & $a,b,c,d$ & $a^*,b^*,c^*,d^*$ \\ \hline
4-4 & $0664287241436146$; & $6,0,-3,9$ & $6,0,-9,3$ \\
& $1177658653747250$ &&\\
4-5 & $0664452175768367$; & $-2,4,5,9$ & $6,0,-9,3$ \\
& $1175561631427380$ &&\\
4-6 & $0664271564363774$; & $0,6,3,9$ & $-8,2,-3,7$ \\
& $1176735233364512$ &&\\
4-7 & $0664463272861135$; & $6,0,9,3$ & $-6,-8,-5,1$ \\
& $1176513423426451$ &&\\
4-8 & $0664151272416748$; & $6,4,-7,5$ & $6,0,-9,3$ \\
& $1176758485764363$ &&\\
5-5 & $0653172153254877$; & $2,4,9,-5$ & $2,4,-9,5$ \\
& $1186526273422531$ &&\\
5-6 & $0664463478185727$; & $-4,2,5,9$ & $-8,2,3,7$ \\
& $1177653314631253$ &&\\
5-7 & $0653487153721511$; & $6,8,-1,5$ & $2,4,5,-9$ \\
& $1186322757876711$ &&\\
5-8 & $0664161572851367$; & $4,6,5,7$ & $4,2,-5,9$ \\
& $1177726542361461$ &&\\
6-6 & $0653151463787817$; & $-2,8,-3,7$ & $2,8,3,-7$ \\
& $1187652534475472$ &&\\
6-7 & $0664277581637113$; & $2,8,3,7$ & $6,8,5,1$ \\
& $1177664515243633$ &&\\
6-8 & $0664475465821113$; & $6,4,5,7$ & $2,8,7,-3$ \\
& $1177546132182722$ &&\\
7-7 & $0664172363751142$; & $8,6,5,-1$ & $8,6,-1,5$ \\
& $1176525365342831$ &&\\
7-8 & $0664475185416311$; & $6,8,1,5$ & $6,4,-5,7$ \\
& $1175416238735363$ &&\\
8-8 & $0664453177857275$; & $-4,6,5,7$ & $4,-6,7,5$ \\
& $1176553815132731$ &&\\
\hline
\end{tabular} \\
\end{center}

\begin{center}
\begin{tabular}{|c|l|l|l|}
\multicolumn{4}{c}{Table 6: $BS(33,32)$} \\ \hline 
Edge & $A$ \& $B$; $C$ \& $D$ & $a,b,c,d$ & $a^*,b^*,c^*,d^*$ \\ \hline
1-1 & $07643661131422181$; & $11,3,0,0$ & $11,3,0,0$ \\
  & $1286331583848171$ &&\\
1-2 & $07644347541711811$; & $3,11,0,0$ & $7,-9,0,0$ \\
  & $1287716551833826$ &&\\
1-3 & $07642434354781830$; & $-1,1,8,8$ & $3,-11,0,0$ \\
  & $1286715346831111$ &&\\
1-4 & $07642414351367712$; & $3,11,0,0$ & $7,-1,8,4$ \\
  & $1284656553724755$ &&\\
1-5 & $07643151228512711$; & $11,3,0,0$ & $-5,-5,8,4$ \\
  & $1287676581466462$ &&\\
1-6 & $07641116654178182$; & $3,11,0,0$ & $7,7,4,4$ \\
  & $1283857157633244$ &&\\
2-2 & $07841512343414140$; & $9,7,0,0$ & $9,7,0,0$ \\
  & $1663752642548557$ &&\\
2-3 & $06613883181363680$; & $1,-1,8,8$ & $9,7,0,0$ \\
  & $1166661118633681$ &&\\
2-4 & $07641411467215623$; & $9,7,0,0$ & $1,7,4,8$ \\
  & $1287676534628461$ &&\\
2-5 & $07644776741834562$; & $-7,9,0,0$ & $5,5,4,8$ \\
  & $1283561165748383$ &&\\
2-6 & $07642431513713560$; & $7,9,0,0$ & $7,-7,4,4$ \\
  & $1283556571663853$ &&\\
3-3 & $07237773326362331$; & $1,1,8,8$ & $1,1,8,8$ \\
  & $1863661181633311$ &&\\
3-4 & $07641562387182580$; & $1,-1,8,8$ & $1,7,8,4$ \\
  & $1285614117616664$ &&\\
3-5 & $07644814118241362$; & $5,5,8,-4$ & $1,1,8,8$ \\
  & $1284626522431467$ &&\\
3-6 & $07786885528463431$; & $-7,-7,4,4$ & $1,1,-8,8$ \\
  & $1286525731546371$ &&\\
4-4 & $07632712148552560$; & $7,1,8,4$ & $7,1,8,4$ \\
  & $1283745543432111$ &&\\
4-5 & $07643457175562810$; & $1,7,8,4$ & $5,-5,4,8$ \\
  & $1283871353112172$ &&\\
4-6 & $07644123143216771$; & $7,7,4,4$ & $7,-1,8,4$ \\
  & $1287661715652463$ &&\\
5-5 & $07641561751648621$; & $5,5,8,4$ & $5,5,8,4$ \\
  & $1285616124737125$ &&\\
5-6 & $07642753664476473$; & $-5,5,8,4$ & $7,-7,4,4$ \\
  & $1285131344465413$ &&\\
6-6 & $07842423683125320$; & $7,-7,-4,-4$ & $7,-7,4,4$ \\
& $1288686675821473$ &&\\
\hline
\end{tabular} \\
\end{center}

\begin{center}
\begin{tabular}{|c|l|l|l|}
\multicolumn{4}{c}{Table 6 (continued)} \\ \hline 
Edge & $A$ \& $B$; $C$ \& $D$ & $a,b,c,d$ & $a^*,b^*,c^*,d^*$ \\ \hline
7-7 & $07632833612216140$; & $11,1,2,-2$ & $11,1,2,2$ \\
  & $1285664844541762$ &&\\
7-8 & $07632578175158550$; & $-1,5,10,2$ & $11,1,2,2$ \\
  & $1283617225865111$ &&\\
7-9 & $07644318776114630$; & $1,11,2,-2$ & $9,3,2,6$ \\
  & $1284842363371533$ &&\\
7-10 & $07645728186111662$; & $3,7,6,6$ & $11,-1,2,2$ \\
  & $1281566243114774$ &&\\
8-8 & $07641431563668731$; & $1,5,2,10$ & $5,1,10,2$ \\
  & $1287676571651331$ &&\\
8-9 & $07632612542858710$; & $5,-1,10,2$ & $9,3,6,2$ \\
  & $1285326563571112$ &&\\
8-10 & $07786157654765620$; & $-3,7,6,6$ & $5,-1,10,2$ \\
  & $1287671165413323$ &&\\
9-9 & $07643428324116160$; & $9,3,2,6$ & $9,3,6,2$ \\
  & $1287761355637215$ &&\\
9-10 & $07644764313231670$; & $3,9,6,2$ & $3,-7,6,6$ \\
  & $1282876155416351$ &&\\
10-10 & $07786231134327142$; & $3,7,6,6$ & $7,3,6,6$ \\
  & $1287335713121563$ &&\\
\hline
\end{tabular} \\
\end{center}

\begin{center}
\begin{tabular}{|c|l|l|l|}
\multicolumn{4}{c}{Table 7: $BS(34,33)$} \\ \hline 
Edge & $A$ \& $B$; $C$ \& $D$ & $a,b,c,d$ & $a^*,b^*,c^*,d^*$ \\ \hline
1-1 & $01643272281847733$; & $2,0,11,3$ & $0,2,3,11$ \\
& $64437112182612640$ &&\\
1-2 & $06426183724377472$; & $0,2,3,11$ & $10,0,-5,3$ \\
& $16715585714616133$ &&\\
1-3 & $01714352388163846$; & $2,0,11,3$ & $4,10,3,3$ \\
& $64462212371615313$ &&\\
1-4 & $06444714358667236$; & $0,2,9,7$ & $2,0,-3,11$ \\
& $16771235115272541$ &&\\
1-5 & $01644816586568712$; & $2,0,-3,11$ & $4,6,9,-1$ \\
& $64715715371825472$ &&\\
1-6 & $01716725382367832$; & $2,0,3,11$ & $8,6,-5,3$ \\
& $64164177317682140$ &&\\
1-7 & $01235326158287371$; & $6,0,7,7$ & $0,2,3,11$ \\
& $66571275124148160$ &&\\
2-2 & $06437211421686264$; & $10,0,3,5$ & $0,-10,3,5$ \\
& $16748465174364121$ &&\\
2-3 & $01846171746522125$; & $10,4,3,-3$ & $0,10,3,5$ \\
& $64311276482826751$ &&\\
2-4 & $01644135625178262$; & $10,0,-5,-3$ & $0,2,7,9$ \\
& $64186575647548281$ &&\\
2-5 & $06552354172663716$; & $6,4,-1,9$ & $0,10,3,5$ \\
& $11746533381536372$ &&\\
2-6 & $06175464158337517$; & $0,10,5,3$ & $-6,8,5,3$ \\
& $12653715652536332$ &&\\
2-7 & $01745856378115355$; & $0,6,7,7$ & $-10,0,3,-5$ \\
& $64775443215125130$ &&\\
3-3 & $06632247171213645$; & $10,4,3,-3$ & $4,10,3,-3$ \\
& $11382325857554252$ &&\\
3-4 & $01738165165617653$; & $4,10,3,-3$ & $2,0,7,9$ \\
& $64433237582218162$ &&\\
3-5 & $06176424834163271$; & $6,4,9,-1$ & $4,10,-3,3$ \\
& $12441571842462262$ &&\\
3-6 & $06441362614513772$; & $8,6,3,5$ & $10,4,3,-3$ \\
& $16776763822163151$ &&\\
3-7 & $06482612536431236$; & $10,-4,3,3$ & $0,6,7,7$ \\
& $12461662575778260$ &&\\
4-4 & $01848235737566316$; & $0,2,7,-9$ & $2,0,7,-9$ \\
& $61242628662324763$ &&\\
4-5 & $06864765526373544$; & $-2,0,9,7$ & $4,6,9,-1$ \\
& $11471612568726141$ &&\\
4-6 & $01644614754247125$; & $8,6,3,5$ & $2,0,7,9$ \\
& $64187352131157381$ &&\\
\hline
\end{tabular} \\
\end{center}

\begin{center}
\begin{tabular}{|c|l|l|l|}
\multicolumn{4}{c}{Table 7 (continued)} \\ \hline 
Edge & $A$ \& $B$; $C$ \& $D$ & $a,b,c,d$ & $a^*,b^*,c^*,d^*$ \\ \hline
4-7 & $06862467734722615$; & $2,0,7,9$ & $0,6,7,-7$ \\
& $11675136251536272$ &&\\
5-5 & $06178545552317721$; & $4,6,1,9$ & $6,-4,1,9$ \\
& $16157375762546143$ &&\\
5-6 & $01837321432341743$; & $6,4,-1,9$ & $8,6,3,5$ \\
& $64381477511564642$ &&\\
5-7 & $01782525345315536$; & $6,0,7,7$ & $4,6,-1,-9$ \\
& $61216753588111360$ &&\\
6-6 & $06554252236661836$; & $8,-6,3,5$ & $6,8,3,5$ \\
& $11257264681477341$ &&\\
6-7 & $01176167385241254$; & $8,6,3,-5$ & $6,0,7,7$ \\
& $66482625745862150$ &&\\
7-7 & $01653673337281734$; & $0,6,-7,7$ & $6,0,-7,7$ \\
& $61473278766448712$ &&\\
\hline
\end{tabular} \\
\end{center}

\begin{center}
\begin{tabular}{|c|l|l|l|}
\multicolumn{4}{c}{Table 8: $BS(35,34)$} \\ \hline 
Edge & $A$ \& $B$; $C$ \& $D$ & $a,b,c,d$ & $a^*,b^*,c^*,d^*$ \\ \hline
1-6 & $076761387537518140$; & $-3,11,-2,-2$ & $-1,-11,0,4$ \\
& $12564355586883173$ &&\\
1-7 & $076544121368256783$; & $1,-1,10,6$ & $11,1,0,-4$ \\
& $12564782516511413$ &&\\
1-8 & $076813248638665463$; & $-3,-5,-2,10$ & $11,1,4,0$ \\
& $12447117137868357$ &&\\
1-9 & $076518224314621143$; & $11,1,0,4$ & $9,7,2,2$ \\
& $12526647176348573$ &&\\
1-10 & $076813687885775451$; & $-11,1,4,0$ & $7,7,6,2$ \\
& $12534322471376448$ &&\\
2-6 & $076813155217615621$; & $9,5,4,4$ & $11,3,-2,2$ \\
& $12534846527651176$ &&\\
2-7 & $076813446761268643$; & $-1,1,10,-6$ & $9,-5,4,4$ \\
& $12534882516521623$ &&\\
2-8 & $076823855753834630$; & $-5,-3,10,-2$ & $5,-9,4,4$ \\
& $12441164836225723$ &&\\
2-9 & $076541256821114362$; & $9,5,4,-4$ & $7,-9,-2,2$ \\
& $12663152258827675$ &&\\
2-10 & $076814553215115552$; & $7,7,2,6$ & $9,5,4,4$ \\
& $12534517187266537$ &&\\
3-6 & $076423483282237882$; & $-3,-11,2,2$ & $3,-1,8,8$ \\
& $12837164638247415$ &&\\
3-7 & $053765656464871261$; & $1,1,-6,10$ & $3,-1,-8,-8$ \\
& $17765746348615187$ &&\\
3-8 & $076544215376333280$; & $3,1,8,-8$ & $-3,-5,10,-2$ \\
& $12662553656248264$ &&\\
3-9 & $076541326141144653$; & $9,7,2,-2$ & $-1,-3,-8,8$ \\
& $12565462867178642$ &&\\
3-10 & $076541313864244753$; & $1,3,8,-8$ & $7,-7,2,6$ \\
& $12565532682263655$ &&\\
4-6 & $076821154786531510$; & $5,7,8,0$ & $11,-3,-2,2$ \\
& $12441254686615465$ &&\\
4-7 & $076542388881587133$; & $-7,-5,0,8$ & $-1,1,10,-6$ \\
& $12634625571747754$ &&\\
4-8 & $076821421676434513$; & $5,3,2,10$ & $7,5,8,0$ \\
& $12441325771765766$ &&\\
4-9 & $076535878535141762$; & $-5,7,0,8$ & $9,-7,2,-2$ \\
& $17677852174231455$ &&\\
4-10 & $076764325821511142$; & $7,7,6,2$ & $5,-7,-8,0$ \\
& $12563712335271855$ &&\\
5-6 & $076541434617337753$; & $-3,11,2,2$ & $7,-3,-8,-4$ \\
& $12565287475625713$ &&\\
5-7 & $076531753465353411$; & $3,7,8,4$ & $1,1,6,10$ \\
& $12456864615253117$ &&\\
\hline 
\end{tabular} \\
\end{center}

\begin{center}
\begin{tabular}{|c|l|l|l|}
\multicolumn{4}{c}{Table 8 (continued)} \\ \hline 
Edge & $A$ \& $B$; $C$ \& $D$ & $a,b,c,d$ & $a^*,b^*,c^*,d^*$ \\ \hline
5-8 & $076543211437821351$; & $7,3,8,-4$ & $5,-3,10,2$ \\
& $12664184625565624$ &&\\
5-9 & $076542443567112150$; & $9,7,-2,2$ & $3,-7,8,4$ \\
& $12634755737233827$ &&\\
5-10 & $076532871428885871$; & $-7,-7,6,2$ & $7,3,8,4$ \\
& $12455284637661614$ &&\\
\hline 
\end{tabular} \\
\end{center}

\begin{center}
\begin{tabular}{|c|l|l|l|}
\multicolumn{4}{c}{Table 9: $BS(36,35)$} \\ \hline 
Edge & $A$ \& $B$; $C$ \& $D$ & $a,b,c,d$ & $a^*,b^*,c^*,d^*$ \\ \hline
1-1 & $066128524558167276$; & $4,-2,11,1$ & $-4,-2,-11,-1$ \\
& $115512428681612272$ &&\\
1-2 & $061752175573814614$; & $4,10,5,1$ & $-4,2,11,-1$ \\
& $123367131555842723$ &&\\
1-3 & $066224581257478141$; & $6,0,9,5$ & $-2,4,-1,11$ \\
& $114662461732721433$ &&\\
1-4 & $065532351825386471$; & $4,-2,11,1$ & $-4,-6,9,3$ \\
& $116425237255181721$ &&\\
1-5 & $061754416162673578$; & $2,8,5,7$ & $2,4,11,1$ \\
& $123357118657181721$ &&\\
2-2 & $061751252386515416$; & $10,4,5,-1$ & $-10,4,-1,5$ \\
& $123355426257165781$ &&\\
2-3 & $016414317335677244$; & $4,10,5,1$ & $0,6,-5,-9$ \\
& $616264816227573640$ &&\\
2-4 & $016815552241223875$; & $10,-4,-1,5$ & $6,-4,9,3$ \\
& $611817576227536871$ &&\\
2-5 & $016622578364127651$; & $8,2,7,5$ & $-4,10,1,-5$ \\
& $615115726753232751$ &&\\
3-3 & $016735472553122818$; & $6,0,5,9$ & $6,0,-9,-5$ \\
& $617256517744612640$ &&\\
3-4 & $065532881558483613$; & $0,-6,5,9$ & $-4,-6,3,-9$ \\
& $116421432717756380$ &&\\
3-5 & $064713642432155468$; & $6,0,9,5$ & $-2,8,-5,7$ \\
& $123473144616418230$ &&\\
4-4 & $068753558343827566$; & $-6,-4,3,9$ & $6,4,9,3$ \\
& $116258162734771362$ &&\\
4-5 & $063828821524645187$; & $2,-8,-7,-5$ & $-6,4,3,9$ \\
& $128664764835184882$ &&\\
5-5 & $016567812318227135$; & $8,2,7,5$ & $8,2,5,7$ \\
& $611513566817626551$ &&\\
\hline 
\end{tabular} \\
\end{center}

\begin{center}
\begin{tabular}{|c|l|l|l|}
\multicolumn{4}{c}{Table 10: $BS(37,36)$} \\ \hline 
Edge & $A$ \& $B$; $C$ \& $D$ & $a,b,c,d$ & $a^*,b^*,c^*,d^*$ \\ \hline
1-1 & $0642483723773112832$; & $1,1,12,0$ & $1,1,-12,0$ \\ 
& $162444213616245723$ &&\\
1-4 & $0781647583615324282$; & $-1,-1,0,12$ & $-1,-9,0,-8$ \\ 
& $167557211545523777$ &&\\
1-5 & $0876588628114455150$; & $1,-1,12,0$ & $9,7,-4,0$ \\
& $161427612285272431$ &&\\
1-6 & $0764841234846532153$; & $3,-3,8,8$ & $-1,1,-12,0$ \\
& $165154775335162126$ &&\\
2-5 & $0767144683434761771$; & $-5,11,0,0$ & $7,-9,-4,0$ \\
& $124873577128343623$ &&\\
2-6 & $0785618342468563210$; & $3,-3,8,8$ & $-5,-11,0,0$ \\
& $126551157157241538$ &&\\
3-4 & $0764214143622153442$; & $11,3,0,4$ & $-1,-9,8,0$ \\
& $164323881543744174$ &&\\
3-5 & $0616123851727712413$; & $9,7,-4,0$ & $-3,11,-4,0$ \\ 
& $123473518825755738$ &&\\ 
3-6 & $0785641356385516141$; & $3,3,8,8$ & $11,3,0,-4$ \\
& $126547124474373121$ &&\\
4-5 & $0616123816854524630$; & $9,-1,8,0$ & $9,7,0,4$ \\
& $123473466224661443$ &&\\
4-6 & $0865126576744588551$; & $-3,-3,8,8$ & $1,9,0,-8$ \\
& $161271417534556246$ &&\\
5-5 & $0717855753413764382$; & $-7,9,0,4$ & $9,-7,0,-4$ \\
& $186871154644661856$ &&\\
5-6 & $0615512388414537671$; & $3,3,8,-8$ & $7,-9,-4,0$ \\ 
& $126528535625368412$ &&\\
7-7 & $0864743671415823362$; & $-1,3,10,-6$ & $-1,3,6,-10$ \\
& $163242244661482565$ &&\\
7-8 & $0778285253655118732$; & $-3,1,10,-6$ & $5,9,2,6$ \\
& $128814612256532345$ &&\\
7-9 & $0764367614152248343$; & $3,1,6,-10$ & $7,5,6,-6$ \\
& $162525438825532618$ &&\\
8-8 & $0762165645151374421$; & $9,5,6,2$ & $5,9,2,-6$ \\ 
& $162127434137824455$ &&\\
8-9 & $0868124566444233641$; & $5,-7,6,6$ & $5,9,2,-6$ \\ 
& $161842144127757326$ &&\\
9-9 & $0637461414423752660$; & $7,5,6,6$ & $7,5,-6,-6$ \\ 
& $128647367258131611$ &&\\
\hline 
\end{tabular} \\
\end{center}

\end{document}